\documentclass[a4paper,12pt]{article}
\usepackage{amsmath,amssymb} \usepackage{euler,eucal}
\usepackage[utf8]{inputenc}
\usepackage[T1]{fontenc}
\usepackage[english]{babel}
\usepackage{tgpagella}
\usepackage[unicode=true,plainpages=false]{hyperref}
\usepackage{graphicx}
\hypersetup{colorlinks=false,pdfborder={0 0 0.1},linkcolor=magenta,anchorcolor=magenta,urlcolor=blue,citecolor=blue,
          pdftitle={One-site density matrix renormalization group and alternating minimum energy algorithm},
          pdfauthor={Sergey V. Dolgov, Dmitry V. Savostyanov},
          pdfkeywords={high-dimensional problems, DMRG, MPS, tensor train format, extreme eigenstate}
          }

\usepackage[left=30mm,right=30mm,top=20mm,bottom=30mm,a4paper]{geometry}
\usepackage{tikz} \usetikzlibrary{positioning,shapes}
\usepackage{pgfplots} \usetikzlibrary{plotmarks}
\usetikzlibrary{calc,3d,shadows}

\def\eps{\varepsilon}
\def\O{\mathcal{O}}
\def\P{\mathcal{P}}
\def\X{\mathcal{X}}

\def\H{\mathbf{H}}
\def\new{\star}
\def\trans{*}
\def\Span{\mathop{\mathrm{span}}\nolimits}

\begin{document}
\author{
  Sergey V. Dolgov\footnotemark[1] ~and
  Dmitry V. Savostyanov\footnotemark[2]
}
\title{One-site density matrix renormalization group and alternating minimum energy algorithm}
\date{December 23, 2013}
\maketitle
\renewcommand{\thefootnote}{\fnsymbol{footnote}}
\footnotetext[1]{Max-Planck Institute for Mathematics in the Sciences, Inselstra{\ss}e 22, Leipzig 04103, Germany  ({\tt sergey.v.dolgov@gmail.com})}
\footnotetext[2]{University of Southampton, School of Chemistry, Highfield Campus, Southampton SO17 1BJ, United Kingdom  ({\tt dmitry.savostyanov@gmail.com})}
\renewcommand{\thefootnote}{\arabic{footnote}}

\begin{abstract}
Given in the title are two algorithms to compute the extreme eigenstate of a high-dimensional Hermitian matrix using the tensor train (TT) / matrix product states (MPS) representation.
Both methods empower the traditional alternating direction scheme with the auxiliary (e.g. gradient) information, which substantially improves the convergence in many difficult cases.
Being conceptually close, these methods have different derivation, implementation, theoretical and practical properties.
We emphasize the differences, and reproduce the numerical example to compare the performance of two algorithms.
\end{abstract}

\section{Introduction} \label{sec:intro}   
Actual problems of science, engineering and society can be so complex, that their mathematical portrait requires more than three dimensions.
Quantum world gives us a perfect example of essentially high--dimensional systems, described by a joint wavefunction (or density matrix) of all particles.
A simple system of $d$ spin-$\frac12$ particles is an entanglement of $\O(2^d)$ possible states, and should be described by the same amount of numbers, which creates out-of-memory errors on a typical workstation for $d \gtrsim 30.$
Even with a brute force of modern supercomputers, standard numerical methods can not honestly simulate protein-size molecules ($d \sim 10^3$ --- $10^4$), since the complexity and storage explode exponentially with $d.$

To overcome this problem, known as the~\emph{curse of dimensionality}, we use data-sparse representations for high-dimensional vectors and matrices, and develop special algorithms to work with them.
Proposed in 1992, the \emph{density matrix renormalization group} (DMRG) algorithm~\cite{white-dmrg-1992} and the \emph{matrix product states} (MPS) formalism~\cite{fannes-mps-1992} suggest to represent a wavefunction $x$ in the following tensor-product form
 \begin{equation}\label{eq:tt}
  \begin{split}
   x  = \tau(x^{(1)},\ldots,x^{(d)})
    & = \sum_{\alpha_1=1}^{r_1} \cdots \sum_{\alpha_{d-1}=1}^{r_{d-1}}
       x^{(1)}_{\alpha_1} \otimes x^{(2)}_{\alpha_1\alpha_2} \otimes \ldots \otimes x^{(d)}_{\alpha_{d-1}}, \\
   x(i_1,\ldots,i_d)
    & = \sum_{\alpha_1=1}^{r_1} \cdots \sum_{\alpha_{d-1}=1}^{r_{d-1}}
       x^{(1)}_{\alpha_1}(i_1) x^{(2)}_{\alpha_1\alpha_2}(i_2) \ldots x^{(d)}_{\alpha_{d-1}}(i_d).
  \end{split}
 \end{equation}
In numerical linear algebra this format was re-discovered as the~\emph{tensor train} (TT) decomposition~\cite{ot-tt-2009,osel-tt-2011}.
A single TT core (or \emph{site}) $x^{(k)}=[x^{(k)}_{\alpha_{k-1}\alpha_k}(i_k)]$ is described by $r_{k-1}n_kr_k$ numbers, where $n_k$ denotes the number of possible states for the $k$--th particle (the \emph{mode size}), and $r_k$ is the TT rank
(or \emph{bond dimension}).
The total number of representation parameters scales as $\O(dnr^2),$ $n\sim n_k,$ $r\sim r_k,$ and is feasible for computations with $d,n,r\lesssim 10^3.$

\begin{figure}[t]
 \colorlet{darkgreen}{green!65!black}
  \hbox to \textwidth{
  \hfill
 \resizebox{.37\textwidth}{!}{\begin{tikzpicture}[x={(5mm,3mm)}, y={(-7mm,1.5mm)}, z={(0mm,10mm)},
                   >=latex,draw opacity=.75,fill opacity=.25]
\tikzset{facestyle/.style={fill=black!10,
                           line join=round, rounded corners=.5mm,thin,double}}

\draw[facestyle,draw=darkgreen] (0,-2,0)  -- (0,2,0);
\draw[facestyle,draw=blue]      (0,0,-.5) -- (0,0,2);
\draw[facestyle,draw=red]      (-1,0,1.5) -- (3,0,1.5);
  
\draw[very thick,darkgreen,->,opacity=.7] (0,1.5,0) -- (0,0,0);
\draw[very thick,blue,->,opacity=.7]      (0,0,0)   -- (0,0,1.5);
\draw[very thick,red,->,opacity=.7]       (0,0,1.5) -- (1.5,0,1.5);

\shadedraw[ball color=darkgreen!60!black,draw opacity=0,fill opacity=.9,ultra thin] ( 0,1.5,0 ) circle (.7mm);
\shadedraw[ball color=darkgreen!60!blue, draw opacity=0,fill opacity=.9,ultra thin] ( 0,0,0 ) circle (.7mm);
\shadedraw[ball color=blue!60!red,   draw opacity=0,fill opacity=.9,ultra thin] (0,0,1.5) circle (.7mm);
\shadedraw[ball color=red!60!black,  draw opacity=0,fill opacity=.9,ultra thin] (1.5,0,1.5) circle (.7mm);
\end{tikzpicture}}
 \hfill
 \resizebox{.37\textwidth}{!}{\begin{tikzpicture}[x={(5mm,3mm)}, y={(-7mm,1.5mm)}, z={(0mm,10mm)},
                   >=latex,draw opacity=.75,fill opacity=.25]
\tikzset{facestyle/.style={fill=black!10,
                           line join=round, rounded corners=.5mm,thin,double}}

\begin{scope}[canvas is yx plane at z=0] 
  \path[facestyle,draw=darkgreen,fill opacity=0] (0,0) rectangle (3,3);
  \path[facestyle,draw=darkgreen] (0,2) rectangle (3,3);
\end{scope}

\begin{scope}[canvas is zx plane at y=1] 
  \path[facestyle,draw=red,fill opacity=0] (1,1) rectangle (2,4);
  \path[facestyle,draw=red] (1,2) rectangle ( 2,4);
  \draw[very thick,red,->,opacity=.9] (1.5,2) .. controls (1.5,3) .. (1.8,3.5);
\end{scope}

\begin{scope}[canvas is zy plane at x=2] 
  \path[facestyle,draw=blue,fill opacity=0] (-1,0) rectangle (2,3);
  \path[facestyle,draw=blue,shade] (0,0) rectangle (2,3);
  \path[facestyle,draw=blue]       (0,0) rectangle (-1,3);
  \draw[very thick,blue,->,opacity=.9] (0,2) .. controls (1,1.8) .. (1.5,1);
\end{scope}

\begin{scope}[canvas is zx plane at y=1] 
  \path[facestyle,draw=red,shade] (1,1) rectangle ( 2,2);
\end{scope}

\begin{scope}[canvas is yx plane at z=0] 
  \path[facestyle,draw=darkgreen,shade] (0,0) rectangle (3,2);
  \draw[very thick,darkgreen,->,opacity=.9] (2.5,.5) .. controls (1.8,1) .. (2,2);
\end{scope}

\shadedraw[ball color=darkgreen!60!black,draw opacity=0,fill opacity=.9,ultra thin] (.5,2.5,0  ) circle (.7mm);
\shadedraw[ball color=darkgreen!60!blue, draw opacity=0,fill opacity=.9,ultra thin] (2 ,2  ,0  ) circle (.7mm);
\shadedraw[ball color=blue!60!red,   draw opacity=0,fill opacity=.9,ultra thin] (2 ,1  ,1.5) circle (.7mm);
\shadedraw[ball color=red!60!black,  draw opacity=0,fill opacity=.9,ultra thin] (3.5,1 ,1.8) circle (.7mm);
\end{tikzpicture}}
  \hfill
  }
\caption{Sequence of low-dimensional optimizations in subspaces $\X_{1},$ $\X_{2},$ $\X_{3},\ldots$ (left), and $\X_{1,2},$ $\X_{2,3},$ $\X_{3,4},\ldots$ (right)}
 \label{fig:3d}
\end{figure}
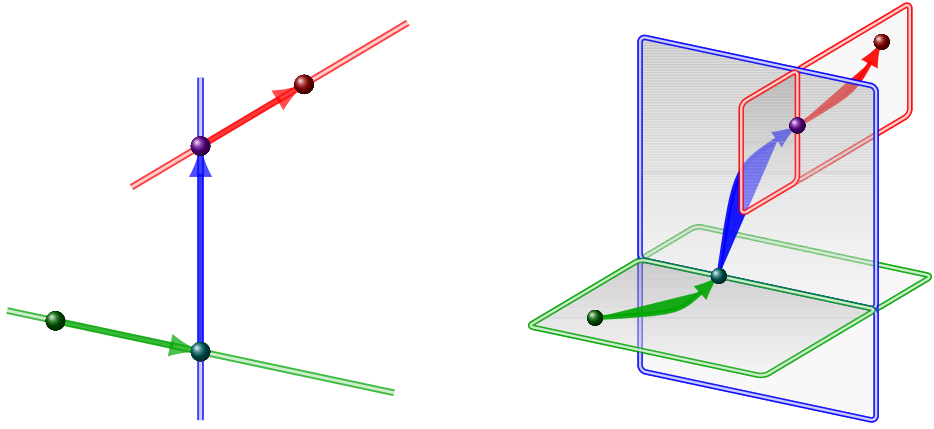

The DMRG algorithm was originally proposed to find the \emph{ground state}, i.e. the minimal eigenpair of a Hermitian matrix $A.$
This problem is equivalent to the minimization of the Rayleigh quotient $J_A(x)=(x,Ax)/(x,x).$
Substituting $J_A(x)$ with $K_A(x)=(x,Ax)-2\Re(x,y),$ and applying the same algorithm, we can solve linear systems $Ax=y$ with Hermitian positive definite matrix~\cite{jeckelmann-dmrgsolve-2002}.
This framework can be extended to a broad class of problems.

Since $x$ is a huge high-dimensional vector, the solution is sought in the structured format~\eqref{eq:tt} with some  TT ranks $r_k,$
defined \emph{a priori} or chosen adaptively.
The simultaneous optimization over all sites is a highly nonlinear and difficult problem.
As it is usual in high-dimensional optimization, we substitute it by a sequence of partial optimizations, each over a particular (small) group of variables.
For our problem, it is natural to group the variables according to the tensor format~\eqref{eq:tt}, e.g. optimize over the components of a  single site $x^{(k)}$ at a time.

The TT format is linear in each site, i.e.
$ 
x = \tau(x^{(1)},\ldots,x^{(d)}) = X_{\neq k} x^{(k)},
$ 
where $X_{\neq k}$ is the $(n_1\ldots n_d) \times (r_{k-1}n_kr_k)$ \emph{frame matrix}, which linearly maps the elements of $x^{(k)}$ to the full vector $x.$
This turns every partial optimization into a local problem of the same type, as the original one,
\begin{equation}\label{eq:loc}
 x^{(k)}_\new = \arg\min_{x^{(k)}} J_A(\tau(x^{(1)},\ldots,x^{(k)},\ldots,x^{(d)}))
              = \arg\min_{x^{(k)}} J_{A_k}(x^{(k)}),
\end{equation}
where $A_k=X_{\neq k}^\trans A X_{\neq k}$ is the $(r_{k-1}n_kr_k)\times(r_{k-1}n_kr_k)$ \emph{reduced matrix}, which inherits the properties of $A,$ i.e. is Hermitian.
Since the frame matrix $X_{\neq k}$ has a structured TT representation (which is the same as~\eqref{eq:tt} with $x^{(k)}$ substituted by the identity matrix), the reduced matrix $A_k$ can be assembled avoiding the exponential costs.
Finally, introducing simple orthogonality conditions for all sites but $x^{(k)},$ we can make the whole matrix $X_{\neq k}$ orthogonal~\cite{schollwock-2011}.
As a consequence, $A_k$ becomes better conditioned than $A,$ and the reduced functional writes
$J_{A_k}(x) = (x,A_k x)/(x,x)$ for the ground state problem, and
$K_{A_k}(x) = (x,A_k x)- 2\Re(x,X_{\neq k}^\trans y)$ for the linear system.
Each optimization~\eqref{eq:loc} is now a classical problem of a tractable size, that can be solved by classical algorithms of numerical linear algebra.

Each local step~\eqref{eq:loc} finds $\min_{x\in\X_{k}}J_A(x),$ where the subspace $\X_{k} = \Span X_{\neq k}$ is of dimension $r_{k-1}n_kr_k$, see Fig \ref{fig:3d} (left).
Here and later by $\Span X$ we denote the subspace of columns of a matrix $X.$
If TT ranks are fixed, the local convergence of such scheme can be analysed using standard methods of multivariate analysis~\cite{ushmaev-tt-2013}.
However, in numerical practice the tensor ranks of the solution are not known in advance, and fixed-rank optimization with wrong ranks would not be efficient.
The DMRG scheme with variable TT ranks is more advantageous, but  the theoretical analysis is even more difficult.

When we allow TT ranks to grow, the dimensions of subspaces $\X_k$ grow as well, and we can use different strategies to expand the subspaces.
Originally, the one-site DMRG scheme (DMRG1) increased the rank $r_k$ by adding (random) orthogonal vectors to $\X_k$, but this algorithm often got stuck far from the ground state.
The problem was solved using \emph{two sites} instead of one in the optimization step~\cite{white-dmrg-1992}.
The \emph{two-site} DMRG algorithm (DMRG2) merges blocks $x^{(k)}$ and $x^{(k+1)},$ and solves the local optimization problem in $\X_{k,k+1} = \Span X_{\neq\{k,k+1\}},$ see Fig.~\ref{fig:3d} (right).
Here $X_{\neq {k,k+1}}$ is the $(n_1\ldots n_d) \times (r_{k-1}n_kn_{k+1}r_{k+1})$  matrix, which has the same TT representation as~\eqref{eq:tt} with blocks $x^{(k)}$ and $x^{(k+1)}$ replaced by the identities.
The DMRG2 converges remarkably well (and is in fact a method of choice) for 1D systems with short--range interactions, but  the cost is approximately $n$ times larger than in the DMRG1.
For systems with long-range interactions two neighboring sites do not provide sufficient information, and DMRG2 can stagnate as well.
To simulate such systems faster and more accurately, better methods to choose search subspaces are required.

The \emph{gradient direction} is central in the theory of optimization methods, and many algorithms use the gradient or its approximate surrogates.
In 2005, S. White proposed the \emph{corrected one-site} DMRG algorithm (DMRG1c),  which adds auxiliary direction to improve the convergence and reduce the computational cost, see \cite{white-dmrg-onesite-2005} and \cite[Sec 6.3]{schollwock-2011} for more details.
In this paper we compare the DMRG1c with the \emph{alternating minimum energy} (AMEn) algorithm.
The AMEn algorithm was recently proposed in~\cite{ds-amr1-2013,ds-amr2-2013} for the solution of linear equations, and the version for the ground state problem appears immediately when we choose $J_A(x)$ as a target function.
In the next section we compare the ideas and implementation aspects of both methods and explain the motivation behind AMEn from numerical linear algebra perspective.
In Sec.~\ref{sec:num} we reproduce a numerical experiment of S. White from~\cite{white-dmrg-onesite-2005}, and demonstrate that AMEn can solve it better than DMRG1c.

\section{Comparison of methods} \label{sec:cmp}  
Both DMRG1c and AMEn combine the local optimization~\eqref{eq:loc} with the step that injects the auxiliary  information.
Both algorithms are \emph{local}, i.e. modify only one block $x^{(k)}$ at a time (cf. the non-local ``ALS$(t+z)$'' algorithm in \cite{ds-amr1-2013}).
Both methods sequentially cycle over TT blocks ($1,2,\ldots,d,d-1,\ldots$).
In the following we assume that \eqref{eq:loc} was just solved for $x^{(k)},$
and consider the step that corrects $x^{(k)}$ before the optimization passes to the next block $x^{(k+1)}.$
This step does not change the vector $x=\tau(x^{(1)},\ldots,x^{(d)})$ (for AMEn), or perturbs it slightly (for DMRG1c), and therefore has a minor direct effect on $J_A(x).$
However, it inserts additional directions to $\Span X_{\neq k+1},$ that improves the convergence of $J_A(x)$ to its global minimum. 

It is crucial how exactly the block $x^{(k)}$ is modified, and which vectors end up in $\Span X_{\neq k+1}$ after that.
In the following we discuss these details, which constitute the main difference between the DMRG1c and the AMEn.

\subsection{Which vector is targeted: $p=A x$ vs. $z=Ax-J_A(x)x$}
Following the \emph{power iteration} method, the DMRG1c algorithm of S. White targets in addition to the solution $x$ the first Krylov vector $p=Ax.$
The AMEn algorithm uses the gradient direction $z=Ax-J_A(x)x.$
In exact arithmetics this makes no difference, since $\Span\{x,p\} = \Span\{x,z\}.$
In practical computations both $p$ and $z$ are perturbed by
  inevitable machine rounding errors,
  perturbations associated with the use of tensor format, and
  additional errors which appear when a surrogate formula (like~\cite[Eq.~$(14)$]{white-dmrg-onesite-2005}) is used to speed up the computations.
The DMRG1c algorithm is derived from perturbation arguments, valid in the vicinity of the minimum of $J_A(x).$
When $x$ approaches the ground state, the angle between $x$ and $p=Ax$ vanishes, and any perturbation in $Ax$ yields a random new direction.
This creates a certain gap between the theory supporting the DMRG1c, and the practice.

Following the \emph{steepest descent} algorithm, the AMEn uses orthogonal vectors $z$ and $x,$ and $\Span\{x,z\}$ is much more stable to perturbations of $z.$
(In general, the Krylov vectors $\{x, Ax, A^2x,\ldots\}$ form an extremely unstable basis, and orthogonalization is crucial.)
The steepest descent algorithm with $z$ substituted by $\tilde z$ converges as long as $(\tilde z,z)>0.$
For the linear systems this fact is elegantly proven in~\cite{kaas-isd-1987}, and the convergence rate of perturbed method is estimated.
An eigenvalue counterpart follows similarly, and the rate of convergence in $\Span\{x,\tilde z\}$ can be estimated from the spectral range of $A.$
This makes the approach implemented in the AMEn algorithm preferable both theoretically and in practice.

\subsection{What is approximated: subspace $\tilde\P$ vs. vector $\tilde z$}
 The computation of full vectors $p=Ax$ and $z=Ax-J(x)x$ is not possible due to their exponentially large size.
 Since $x$ and $A$ are both in TT format, we can avoid the curse of dimension and represent $p=Ax$ and $z=Ax-J(x)x$ by the TT format.
 However, the TT ranks of $Ax$ can be as large as product of TT ranks of $A$ and $x,$ which makes the calculations difficult.

To reduce these costs, S. White suggests in the {\bf{DMRG1c}} the following scheme.
The TT format~\eqref{eq:tt} is divided in two parts:
 left blocks (number $1,\ldots,k$) are referred to as \emph{system}, and
 right blocks ($k+1,\ldots,d$) as~\emph{environment}.
The TT format for the matrix $A$ is written accordingly,
\begin{equation}\label{eq:a}
 A
   = \sum\limits_{\gamma_1 \ldots \gamma_{d-1}}
   \underbrace{A^{(1)}_{\gamma_1} \otimes \ldots \otimes A^{(k)}_{\gamma_{k-1}\gamma_k}}_{\text{system}}
   \otimes
   \underbrace{A^{(k+1)}_{\gamma_k\gamma_{k+1}} \otimes \ldots \otimes A^{(d)}_{\gamma_{d-1}}}_{\text{environment}},
\end{equation}
or shortly $A = \sum_{\gamma} A^<_{\gamma} \otimes A^>_{\gamma}.$
Similarly, Eq.~\eqref{eq:tt} reduces to $x=\sum_\alpha x^<_\alpha \otimes x^>_\alpha.$
The targeting of $p=Ax$ is substituted by targeting of all $p_\gamma=(A^<_\gamma \otimes I)x.$

Although in general $p\notin\cup_\gamma\Span p_\gamma,$ it can be argued that the set $\{p_\gamma\}$ contains a sufficient subspace information.
To show this, we write
\begin{equation}\label{eq:p}
   p  = \sum\nolimits_{\alpha,\gamma}
                  \left( A^<_\gamma x^<_\alpha \right) \otimes
                  \left( A^>_\gamma x^>_\alpha \right), \qquad
  p_\gamma = \sum\nolimits_{\alpha}
                  \left( A^<_\gamma x^<_\alpha \right) \otimes
                                                          x^>_\alpha,
\end{equation}
and consider vectors $p$ and $p_\gamma$ as system-by-environment matrices $P$ and $P_\gamma$ of size $(n_1\ldots n_k)\times(n_{k+1}\ldots n_d).$
Now $\Span P \subset \cup_\gamma\Span P_\gamma=\P,$ where $A^>_\gamma$ contains the coefficients of the required linear combination --- in the exact arithmetics the \emph{system}-related components of $p$ belong to $\P.$
Each $p_\gamma$ is easier to compute than $p,$ because it does not depend on the environment part $A^>_\gamma.$

The total dimension of $\P$ grows in each step, and to keep TT ranks and storage moderate, we have to truncate it.
The approximation step in the DMRG1c replaces $\P$ with a subspace $\tilde \P$ of a smaller dimension, using a classical \emph{singular value decomposition} (SVD), or Schmidt decomposition technique.
The \emph{dominant subspace} $\tilde\P$ is spanned by the first singular vectors of the matrix $\begin{bmatrix} X & \sqrt{a_1}P_1 & \sqrt{a_2}P_2 & \ldots\end{bmatrix},$
where all target vectors are concatenated with empirically chosen weighting coefficients $a_{\gamma}$.
The method assumes that the vector $p=Ax$ is likely to belong to $\tilde\P.$

This assumption makes perfect sense if $p$ is a random sample from $\P$ --- for a random $u$, $Xu$ is more likely to end up in the dominant subspace of $X.$
However, the target vector $p=Ax$ does not belong to $\tilde\P$ in general, for any choice of weights $\sqrt{a_\gamma}.$
The reason is that $p$ depends crucially on $A^>_\gamma,$ whereas this information is dropped for the sake of faster computations in $p_\gamma$ and hence $\P$ and $\tilde\P$.
Selecting $A^>_\gamma$ in~\eqref{eq:p}, we may come across any vector in $\P$, even the smallest singular vector.
That is, for each choice of $\sqrt{a_\gamma}$ and $x$ there is a \emph{counterexample} of a Hamiltonian, for which the slightest truncation of $\P$ loses the system-related part of the target vector $p=Ax.$

The {\bf{AMEn}} approximates $z=Ax-J_A(x)x$ into its own TT format using any compression tool.
Either the SVD--based technique, which computes the approximation $\tilde z \approx z$ up to any prescribed tolerance $\eps$,
or a faster (but heuristic) \emph{alternating least squares} (ALS) method may be used.
In any case, we may generate an approximation $\tilde z$ with a desired accuracy, which
guarantees the convergence of the steepest descent method with the imperfect direction $\tilde z.$
This fact provides the theoretical bounds for the global convergence rate of the whole AMEn scheme, similarly to~\cite{ds-amr2-2013}.

\subsection{How the new direction is used: averaging vs. enrichment}
The last but not the least detail is how exactly the information about the auxiliary direction is injected in the algorithm.
To show this in isolation from the other dissimilarities outlined above, we assume that in both methods we target in addition to $x$ \emph{only one} vector $s.$
To simplify the presentation we also consider the $d=2$ case, and write
$
x=\sum\nolimits_{\alpha=1}^{r_x} x^<_\alpha x^>_\alpha,
$
and
$
s=\sum\nolimits_{\beta=1}^{r_s} s^<_\beta s^>_\beta,
$
where ``$<$'' and ``$>$'' denote the first and the second blocks, respectively. %

The DMRG1c~\emph{averages} the subspaces
$X=\begin{bmatrix}x^<_1 & \ldots x^<_{r_x}\end{bmatrix}$
and
$S=\begin{bmatrix}s^<_1 & \ldots s^<_{r_s}\end{bmatrix}$
by computing the dominant subspace $\Span U$ of the Gram matrix as follows, $G = X X^\trans + a S S^\trans \approx U U^\trans,$
where
$U=\begin{bmatrix}u^<_1 & \ldots u^<_{r_u}\end{bmatrix}.$
As shown in previous subsection, this procedure does not guarantee that $x$ or $s$ ends up in $\Span (U \otimes I),$ unless $\Span U = \Span\begin{bmatrix} X & S \end{bmatrix}.$
The TT core $x^<$ is replaced by the vectors of $U,$ that introduces a $\O(\sqrt{a})$ perturbation to $x$ and probably worsen $J_A(x).$
It is clear though that $a$ should vanish when we approach the exact solution, but the general recipe is not known.

The AMEn avoids the outlined difficulties by~\emph{merging} $U=\begin{bmatrix} X & S \end{bmatrix}$ and zero-padding the second block.
Values of $x$ and $J_A(x)$ are preserved, no rescaling is required, and both $\{x,s\}\in\Span(U\otimes I).$
The downside is that we choose $r_u=r_x+r_s$ each time we expand the subspaces.
However, when we use the approximate gradient direction $s=\tilde z\approx z=Ax-J_A(x)x$ the low-rank $\tilde z$ usually suffice, e.g. with $r_s\approx r_x/2.$
We can also truncate the TT-ranks at the end of each iteration and control the perturbation to $J_A(x).$

\section{Numerical example} \label{sec:num} 
Following S. White~\cite{white-dmrg-onesite-2005}, we consider the spin-1 periodic Heisenberg chain,
\begin{equation}\label{eq:H}
 \begin{split}
     A  & = \H_1\cdot\H_2 + \H_2\cdot\H_3 + \ldots + \H_{d-1}\cdot\H_d + \H_d\cdot\H_1,                \\
   \H_i & = (\H_i^x, \H_i^y, \H_i^z)^\top, \qquad \H_i\cdot\H_j=H_i^xH_j^x + H_i^yH_j^y + H_i^zH_j^z,  \\
   H^{\{x,y,z\}}_i & = I \otimes\cdots\otimes I \otimes S_{\{x,y,z\}} \otimes I \otimes\cdots\otimes I, \quad\text{$S$ in position $i,$}
 \end{split}
\end{equation}
where  $S_{\{x,y,z\}}$, are the $3 \times 3$ Pauli matrices for spin-$1$ particles.
The number of spins $d$ is set to $100$, i.e. the wavefunction belongs to the $3^{100}$-dimensional Hilbert space.
This example is particularly illustrating, since the mismatch between the linear TT model~\eqref{eq:tt}  and the cycle structure of~\eqref{eq:H} complicates the problem --- the solution has large TT ranks, and both the one-- and two--site DMRG converge slowly.

\begin{figure}[t]
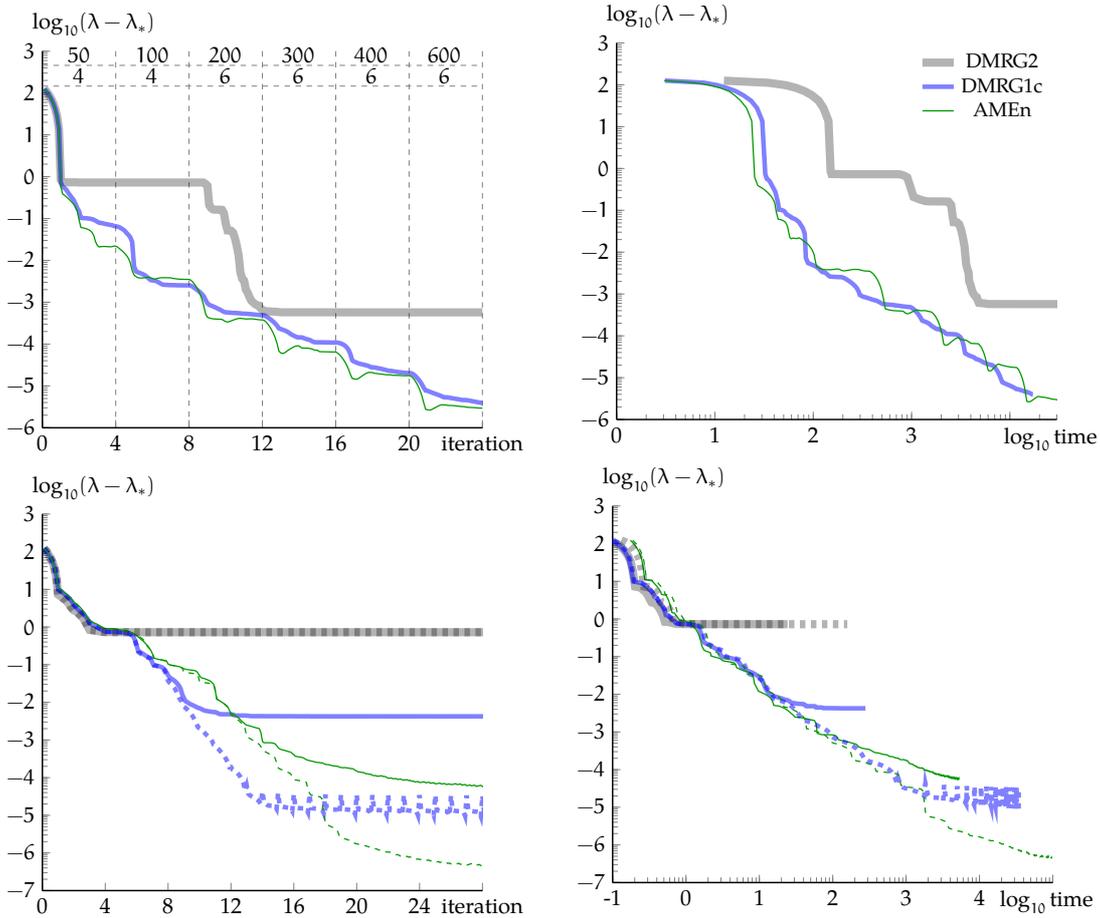

  \resizebox{.48\textwidth}{!}{\input{./pgfart.sty} \begin{tikzpicture}

\pgfplotsset{grid style={dashed,gray}}

\begin{axis}[%
xmode=normal,ymode=log,
xlabel={iteration},
ylabel={$\log_{10}(\lambda-\lambda_{*})$},
xmin=0, xmax=24.01,
xtick={0,4,8,12,16,20,24},
xticklabels={0,4,8,12,16,20,~},
xmajorgrids=true,
ymin=1e-06, ymax=1e3,
legend style={at={(0.99,0.99)},anchor=north east}]

\addplot table [header=false, x index = 0, y index = 6]{./dat/r600.dat};
\addplot table [header=false, x index = 0, y index = 4]{./dat/r600.dat};
\addplot table [header=false, x index = 0, y index = 2]{./dat/r600.dat};
\end{axis}

\path[dashed,draw=gray]       (0.0,7.0) to (8.4,7.0);
\node[black] at (0.7,7.2) {$50$};
\node[black] at (2.1,7.2) {$100$};
\node[black] at (3.5,7.2) {$200$};
\node[black] at (4.9,7.2) {$300$};
\node[black] at (6.3,7.2) {$400$};
\node[black] at (7.7,7.2) {$600$};

\path[dashed,draw=gray]       (0.0,6.6) to (8.4,6.6);
\node[black] at (0.7,6.8) {$4$};
\node[black] at (2.1,6.8) {$4$};
\node[black] at (3.5,6.8) {$6$};
\node[black] at (4.9,6.8) {$6$};
\node[black] at (6.3,6.8) {$6$};
\node[black] at (7.7,6.8) {$6$};
\end{tikzpicture}%
} \hfil
  \resizebox{.48\textwidth}{!}{\input{./pgfart.sty} \begin{tikzpicture}

\begin{axis}[%
xmode=log,ymode=log,
xlabel={$\log_{10}\mathrm{time}$},
ylabel={$\log_{10}(\lambda-\lambda_{*})$},
xmin=1, xmax=3e4,
xtick={1,1e1,1e2,1e3,1e4},
xticklabels={0,1,2,3,~},
ymin=1e-06, ymax=1e3,
legend style={at={(0.99,0.99)},anchor=north east}]

\addplot table [header=false, x index = 5, y index = 6]{./dat/r600.dat};
\addplot table [header=false, x index = 3, y index = 4]{./dat/r600.dat};
\addplot table [header=false, x index = 1, y index = 2]{./dat/r600.dat};

\legend{DMRG2,DMRG1c,AMEn};

\end{axis}
\end{tikzpicture}%
} \\
  \resizebox{.48\textwidth}{!}{\input{./pgfart.sty} \begin{tikzpicture}

\begin{axis}[%
xmode=normal,ymode=log,
xlabel={iteration},
ylabel={$\log_{10}(\lambda-\lambda_{*})$},
xmin=0, xmax=28.01,
xtick={0,4,8,12,16,20,24,28},
xticklabels={0,4,8,12,16,20,24,~},
ymin=1e-7, ymax=1e3,
ytick={1e-7,1e-6,1e-5,1e-4,1e-3,1e-2,1e-1,1e0,1e1,1e2,1e3},
yminorticks=true,
legend style={at={(0.99,0.99)},anchor=north east}]

\addplot table [header=false, x index = 0, y index = 6]{./dat/e3.dat};
\addplot table [header=false, x index = 0, y index = 4]{./dat/e3.dat};
\addplot table [header=false, x index = 0, y index = 2]{./dat/e3.dat};

\addplot+[dashed] table [header=false, x index = 0, y index = 6]{./dat/e4.dat};
\addplot+[dashed] table [header=false, x index = 0, y index = 4]{./dat/e4.dat};
\addplot+[dashed] table [header=false, x index = 0, y index = 2]{./dat/e4.dat};



\end{axis}
\end{tikzpicture}%
} \hfil
  \resizebox{.48\textwidth}{!}{\input{./pgfart.sty} \begin{tikzpicture}

\begin{axis}[%
xmode=log,ymode=log,
xlabel={$\log_{10}\mathrm{time}$},
ylabel={$\log_{10}(\lambda-\lambda_{*})$},
xmin=1e-1, xmax=1e5,
xtick={1e-1,1e0,1e1,1e2,1e3,1e4,1e5},
xticklabels={-1,0,1,2,3,4,~},
ymin=1e-7, ymax=1e3,
ytick={1e-7,1e-6,1e-5,1e-4,1e-3,1e-2,1e-1,1e0,1e1,1e2,1e3},
yminorticks=true,
legend style={at={(.5,-.3)},anchor=north}]

\addplot table [header=false, x index = 5, y index = 6]{./dat/e3.dat};
\addplot table [header=false, x index = 3, y index = 4]{./dat/e3.dat};
\addplot table [header=false, x index = 1, y index = 2]{./dat/e3.dat};

\addplot+[dashed] table [header=false, x index = 5, y index = 6]{./dat/e4.dat};
\addplot+[dashed] table [header=false, x index = 3, y index = 4]{./dat/e4.dat};
\addplot+[dashed] table [header=false, x index = 1, y index = 2]{./dat/e4.dat};


\end{axis}
\end{tikzpicture}%
} \hfil
\caption{Error in the eigenvalue vs. iteration (left) and CPU time (right). Methods: AMEn \cite{ds-amr1-2013,ds-amr2-2013}, DMRG1c \cite{white-dmrg-onesite-2005}, DMRG2 \cite{white-dmrg-1992}. Top: parameters depend on iteration as shown on top of left figure (ranks and $\log_{10}(1/a)$, resp.).
Bottom: $a=10^{-4}$, ranks depend on accuracies: $\eps=10^{-3}$ (solid lines), $\eps=10^{-4}$ (dashed lines)}
\label{fig:heisen}
\end{figure}

The way how the TT ranks are chosen during the algorithm is also very important.
We first adopt the rank selection strategy from \cite{white-dmrg-onesite-2005}, and compare the DMRG2, the DMRG1c and the AMEn algorithms. The results are shown in Fig.~\ref{fig:heisen} (top left), which overlays~\cite[Fig. 3]{white-dmrg-onesite-2005} with the AMEn behavior.
In Fig.~\ref{fig:heisen} (top right) the convergence of $\lambda=J_A(x)$ to the reference value $\lambda_{*}=-140.14840390392$ (computed in \cite{white-dmrg-onesite-2005} by the DMRG1c with TT ranks $4000$) is given w.r.t. the cumulative CPU time.

We see that both DMRG methods correctly reproduce the experiment from \cite{white-dmrg-onesite-2005}: the two-site DMRG stagnates at a high error level, while the corrected DMRG converges significantly faster.
The AMEn method manifests practically the same efficiency.
Since it searches in a larger subspace, it is even more accurate w.r.t. iterations, but becomes slightly slower during the optimization of inner TT blocks.
However, letting it to increase the ranks (each fourth iteration) yields sharper error decays.

To free the algorithm from tuning parameter, we prefer to choose the ranks adaptively to the desired accuracy.
With this we also avoid artificial rank limitation, which pollutes the convergence.
Therefore, in the second experiment we use the same algorithms but perform the truncation of TT blocks via the SVD using the relative Frobenius-norm accuracies $\varepsilon=10^{-3}$ and $\varepsilon=10^{-4}$.
The results are shown in Fig. \ref{fig:heisen} (bottom).

We see that when ranks are chosen adaptively, the AMEn rapidly becomes faster than the other algorithms.
Even the DMRG1c stagnates relatively early, since the correction $p_{\gamma}$ \eqref{eq:p} contaminates the dominant basis of the ground state.
Moreover, since the $\varepsilon$-truncation eliminates the correction if $a \lesssim \varepsilon^2$, it is worthless to decrease the scale $a$ (cf. Fig \ref{fig:heisen}, top left).
Both the adaptivity and speed speak in favour of such truncation: the same accuracy levels are achieved several times faster than in the fixed-rank experiment (e.g. 10 vs. 100 sec. for $\lambda-\lambda_{*}\approx 10^{-2}$ and $\varepsilon=10^{-3}$).
Larger time spent by AMEn in the latter iterations is compensated by a significantly better accuracy, which is close to the optimal level $\O(\varepsilon^2)$.

Finally, the AMEn is applicable to a wider class of problems, and is free from heuristic parameters.

\def\cprime{$'$}

\end{document}